\documentclass[12pt,final]{article}
\usepackage{amsmath,amsthm,amsfonts,amssymb,graphicx,enumerate,psfrag}
\usepackage{mathrsfs}
\usepackage{fullpage}

\newcommand{\cA}{\mathcal {A}}
\newcommand{\cP}{\mathcal {P}}
\newcommand{\dd}{\mathrm {d}}
\newcommand{\dH}{\mathrm {H}\,}
\newcommand{\cF}{\mathcal {F}}
\newcommand{\dX}{\mathbb X}
\newcommand{\cE}{\mathcal {E}}
\newcommand{\dR}{\mathbb {R}}
\newcommand{\tmix}{{\mathrm{t}_{\textsc{mix}}}}
\newcommand{\dtv}{{\mathrm{d}_{\textsc{tv}}}}
\newtheorem{theorem}{Theorem}
\newtheorem{remark}{Remark}
\title{Hyper-contractivity and entropy decay in discrete time}
\author{Justin Salez}
\begin{document}
\maketitle
\abstract{Consider a measure-preserving transition kernel $T$ on an arbitrary probability space $(\dX,\cA,\pi)$. In this level of generality, we prove that a one-step hyper-contractivity estimate of the form $\|T\|_{p\to q}\le 1$ with  $p<  q$  implies a one-step entropy contraction estimate of the form $\dH(\mu T\,|\,\pi)\le \theta\, \dH(\mu\,|\,\pi)$, with $\theta=p/q$. Neither reversibility, nor any sort of regularity is required. This \emph{static} implication is simultaneously simpler and stronger than the  celebrated \emph{dynamic} relation between exponential hyper-contractivity and exponential entropy decay along continuous-time Markov semi-groups. 
}
\section{Introduction}

\paragraph{Transition kernels.} Throughout this note, we consider a measure-preserving transition kernel on a probability space $(\dX,\cA,\pi)$, i.e. a map $T\colon \dX\times\cA\to[0,1]$ such that
\begin{enumerate}[(i)]
\item  $A\mapsto T(x,A)$ is a probability measure for each $x\in \dX$;
\item  $x\mapsto T(x,A)$ is measurable for each $A\in \cA$;
\item $\pi$ is fixed by the natural action of $T$ on $\cP(\dX)$, i.e. the map  $\mu\mapsto\mu T$ given by
\begin{eqnarray}
\label{def:muT}
(\mu T)(A) & := & \int_\dX T(x,A)\,\mu(\dd x).
\end{eqnarray}
\end{enumerate}
Our aim is to shed a new light on the interplay between two fundamental regularization properties of $T$:  \emph{hyper-contractivity} on the one hand, and \emph{entropy contraction} on the other. Let us first briefly recall what those notions are. 

\paragraph{Hyper-contractivity.} The transition kernel $T$ naturally acts on  non-negative measurable functions $f\colon \dX\to[0,\infty]$ via the familiar formula
\begin{eqnarray}
\label{def:T}
 Tf(x) & := & \int_\dX f(y)T(x,\dd y).
\end{eqnarray}
This definition of course extends to signed functions by linearity, as long as $Tf_+$ or $Tf_-$ is finite. Moreover, for each $p\ge 1$, Jensen's inequality and the stationarity property $\pi T=\pi$ easily and classically guarantee that the above action is a contraction on the Banach space $L^p(\pi)$, equipped with its usual norm, 
\begin{eqnarray*}
\|f\|_p & := & \left(\int |f|^p\,\dd\pi\right)^{1/p}.
\end{eqnarray*}
Hyper-contractivity is the stronger requirement that, for some $q>p$ and all $f\in L^p(\mu)$,
\begin{eqnarray}
\label{def:hyper-contractivity}
\|Tf\|_q & \le & \|f\|_p.
\end{eqnarray}
The first estimates of this form were discovered by Nelson along the Ornstein-Uhlenbeck semigroup \cite{MR210416}, and by Bonami and Beckner on the discrete hypercube \cite{MR283496,MR385456}. Following the foundational contributions of Gross \cite{MR420249}, Bakry and Émery \cite{MR889476}, and Diaconis and Saloff-Coste \cite{MR2341319}, hypercontractivity has emerged as a fundamental tool in the quantitative study of Markov processes. In recent years, its impact has extended dramatically, yielding remarkable advances in statistical physics  and computer science \cite{MR4732718,10.1214/24-PS27}. 

\paragraph{Entropy contraction.} The second regularization property that we shall consider is a classical entropy contraction estimate, which takes the form
\begin{eqnarray}
\label{def:entropycontraction}
\forall \mu\in \cP(\dX),\qquad \dH(\mu T\,|\,\pi) & \le & \theta\, \dH(\mu\,|\,\pi),
\end{eqnarray}
for some constant $\theta<1$. Here,  $\dH(\cdot\,|\,\cdot)$ denotes the  Kullback-Leibler divergence:
 \begin{eqnarray*}
\dH(\mu\,|\,\pi) & := & 
\left\{
\begin{array}{ll}\displaystyle{
\int \log \left(\frac{\dd\mu}{\dd\pi}\right)\,\dd \mu} & \textrm{if }\mu\ll\pi\\
+\infty & \textrm{else.}
\end{array}
\right.
 \end{eqnarray*}
Among several other applications, entropy contraction plays a prominent role in the analysis of \emph{mixing times} of Markov processes \cite{MR2283379,MR2341319}, as well as in quantifying  the  celebrated \emph{concentration-of-measure phenomenon} under the reference law $\pi$. We refer the unfamiliar reader to the recent lecture notes \cite{entropy,salez2025modernaspectsmarkovchains} and the references therein  for a self-contained introduction, and many examples.

\section{Result and discussion}

In the present note, we establish the following simple, general, and seemingly new quantitative relation between hyper-contractivity and entropy contraction. 
\begin{theorem}[Main result]\label{th:main}
For any parameters $1\le p\le q$,  the hyper-contractivity estimate (\ref{def:hyper-contractivity}) implies the entropy contraction estimate (\ref{def:entropycontraction}), with $\theta=p/q$.
\end{theorem}

We emphasize that Theorem \ref{th:main} applies to any measure-preserving transition kernel on any probability space: neither reversibility, nor any sort of regularity is required. In particular, we may think of $T$ as the transition kernel of a discrete-time Markov process on $\dX$ with invariant law $\pi$, in which case the conclusion can readily be iterated to provide a geometric rate of convergence to equilibrium, in relative entropy. 
Alternatively, we can choose $T=P_t$, where $(P_t)_{t\ge 0}$ is a given measure-preserving Markov semi-group on $(\dX,\cA,\pi)$ and  $t\ge 0$ a particular time-scale which we want to investigate. In this well-studied continuous-time setting, the ability to focus on a single instant appears to be new, and makes our static implication stronger than its celebrated dynamic counterpart,  which we now review.

To lighten our discussion, we deliberately omit technical details and refer the interested reader to the excellent references \cite{MR3155209} (for Markov diffusions on Euclidean spaces or smooth manifolds) and \cite{MR1410112} (on finite state spaces). Consider  a measure-preserving Markov semi-group   $(P_t)_{t\ge 0}$ on a probability space $(\dX,\cA,\pi)$, and assume that it satisfies an exponential hyper-contractivity estimate of the form
\begin{eqnarray}
\label{assume:beta}
\forall f\in L^2(\pi),\quad\forall t\ge 0,\quad \|P_tf\|_{1+e^{4\beta t}} & \le & \|f\|_2,
\end{eqnarray}
for some  $\beta>0$. Then, a classical differentiation  leads to the \emph{log-Sobolev inequality}
\begin{eqnarray}
\label{LSI}
\forall f\in\cF,\quad \cE(\sqrt{f},\sqrt{f}) & \ge & \beta\, \dH(f\dd \pi\,|\,\pi),
\end{eqnarray}
where $\cF$ is an appropriate class of probability densities on $(\dX,\cA,\pi)$,  and  $\cE(\cdot,\cdot)$ the  Dirichlet form associated with the semi-group. Now, in view of the elementary estimate
\begin{eqnarray*}
\forall a,b\in(0,\infty),\qquad b\log\frac{b}{a} & \ge  & 2\sqrt{b}(\sqrt{b}-\sqrt{a}),
 \end{eqnarray*} the log-Sobolev inequality (\ref{LSI}) always implies its ``modified'' version
\begin{eqnarray}
\label{MLSI}
\forall f\in\cF,\quad \cE(f,\log f)&  \ge & 2\beta\,\dH(f\dd \pi\,|\,\pi),
\end{eqnarray}
which, by a Grönwall-type argument, finally guarantees the exponential entropy decay 
\begin{eqnarray}
\label{assume:alpha}
\forall \mu\in\cP(\dX),\quad\forall t\ge 0,\quad \dH(\mu P_t\,|\, \pi) & \le & e^{-2\beta t}\dH(\mu \,|\, \pi).
\end{eqnarray}
In other words, the implication $(\ref{assume:beta})\Longrightarrow(\ref{assume:alpha})$ always holds along Markov semi-groups. This general relation between hyper-contractivity and entropy contraction is of course a well-established fact, with many applications. It is important to realize, however, that it is inherently dynamical: the estimates $(\ref{assume:beta})$ and $(\ref{assume:alpha})$ hold for all $t\ge 0$, and the  semi-group structure is crucially used to reduce them to the respective functional inequalities (\ref{LSI}) and (\ref{MLSI}), which can then be appropriately compared. In contrast, fixing a particular time $t\ge 0$ and choosing $T=P_t$ in Theorem \ref{th:main} directly yields
\begin{eqnarray*}
\forall \mu\in\cP(\dX),\quad \dH(\mu P_t\,|\, \pi) & \le & \frac{2}{1+e^{4\beta t}}\,\dH(\mu \,|\, \pi),
\end{eqnarray*}
which, in view of the inequality $1+e^{4\beta t}\ge 2e^{\beta t}$,  is always stronger than (\ref{assume:alpha}).  On finite state spaces for example, this readily leads to the mixing-time estimate
\begin{eqnarray*}
\forall\varepsilon\in(0,1),\quad \tmix(\varepsilon) & \le & \frac{1}{4\beta}\left[\log\log\frac{1}{\pi_{\star}}+\log\frac{1}{\varepsilon^2}\right],
\end{eqnarray*}
which is  twice better than what the traditional estimate (\ref{assume:alpha}) would give. Here, we have used the classical notation $\pi_{\star}:=\min_{x\in \dX}\pi(x)$ for the minimum stationary mass, and 
\begin{eqnarray*}
\tmix(\varepsilon) & := & \min\left\{t\ge 0\colon\ \forall\mu\in\cP(\dX),\ \dtv(\mu P_t,\pi)\le\varepsilon\right\},
\end{eqnarray*}
for the worst-case total-variation mixing time. More importantly, our result does not require the assumption \eqref{assume:beta} to hold at all times: a \emph{static} hyper-contractivity estimate at some fixed time suffices to guarantee an entropy contraction estimate at the very same time. Exploiting this instantaneous relation can be quite interesting in practice, because the large-time regularizing effect of $P_t$ is often much stronger than what an infinitesimal computation at $t=0$ would predict. 

Finally,  in addition to being more general and perhaps more natural, Theorem \ref{th:main} admits an elementary proof by duality, which completely avoids the use of semi-groups, hence  the technical precautions needed in order to  safely differentiate along them.

\begin{remark}[Converse] It is well known that the implication  $(\ref{assume:beta})\Longrightarrow(\ref{assume:alpha})$  can  be reversed for Markov diffusions, and an approximate version of this was recently established on discrete state spaces as well \cite{MR4620718,salez2025intrinsicregularitydiscretelogsobolev}. In light of this, it is natural to ask for an appropriate converse to our main theorem. More precisely, if a reversible transition kernel $T$ satisfies a one-step entropy contraction estimate, and if it is sufficiently ``regular'', can one deduce a one-step hyper-contractivity estimate, at a reasonable price? 
\end{remark}

\section{Proof of the theorem}

Let us start by recalling that $T$ has an adjoint $T^\star$, characterized by the duality relation
\begin{eqnarray}
\label{def:adjoint}
\int gT^\star f\,\dd\pi & = & \int fTg \,\dd\pi,
\end{eqnarray}
for any measurable functions $f,g\colon \dX\to\dR$ such that those integrals make sense. Now, fix $1\le p\le q$ and suppose that $T$ satisfies the hyper-contractivity property
\begin{eqnarray}
\label{assume}
\forall g\in L^p(\pi),\qquad \|Tg\|_q & \le & \|g\|_p.
\end{eqnarray}
Given a probability measure $\mu\in\cP(\dX)$, our goal is to prove that
\begin{eqnarray}
\label{goal}
\dH(\mu T\,|\,\pi) & \le & \frac{p}{q}\,\dH(\mu\,|\,\pi).
\end{eqnarray}
We may assume that the right-hand side is finite, otherwise there is nothing to prove. In other words, $\mu$ admits a density $f$ w.r.t. $\pi$, and $f\log f\in L^1(\pi)$. First, using the definition of $\mu T$ at (\ref{def:muT}) and the duality relation (\ref{def:adjoint}) with $g={\bf 1}_A$, we find
\begin{eqnarray*}
\forall A\in\cA,\qquad (\mu T)(A) & = &  \int  f T{\bf 1}_A \,\dd\pi \ = \ \int  {\bf 1}_A T^\star f\,\dd\pi ,
\end{eqnarray*}
which shows that $\mu T$ is also absolutely continuous w.r.t. $\pi$, with density $T^\star f$. Next, the hyper-contractivity assumption (\ref{assume}) applied to $g=(T^\star f)^{1/p}$  reads
\begin{eqnarray*}
\int \left(Te^{\frac{1}{p}\log T^\star f}\right)^q \dd \pi & \le & 1.
\end{eqnarray*}
On the other hand, Jensen's inequality ensures that $e^{Th}\le Te^h$ for any measurable function $h\colon \dX\to\dR$ such that $Th_+<\infty$, hence in particular for $h:=\frac{1}{p}\log T^\star f$. Inserting this pointwise estimate into the above integral, we deduce that the function
\begin{eqnarray*}
\varphi & := & \frac{q}{p} T\log T^\star f,
\end{eqnarray*}
satisfies $\int e^{\varphi}\,\dd \pi \le  1$. Finally, by the variational formulation of entropy (or just the convexity estimate $u\log u\ge 1-u$ applied to $u=fe^{-\varphi}$), this last condition implies
\begin{eqnarray*}
\int f\varphi\, \dd \pi & \le &  \int f\log f\,\dd\pi.
\end{eqnarray*}
The right-hand side is  $\dH(\mu\,|\,\pi)$, and the left-hand side equals $\frac{q}{p}\dH(\mu T\,|\,\pi)$ because
\begin{eqnarray*}
\int f T(\log T^\star f)\,\dd \pi & = & \int (T^\star f)\log (T^\star f)\,\dd \pi \ = \ \dH(\mu T\,|\,\pi),
\end{eqnarray*}
where we have used the duality relation (\ref{def:adjoint}) with $g=\log T^\star f$. Thus, (\ref{goal}) is established.
\section*{Acknowledgment} The author warmly thanks Liming Wu for raising the question answered in the present note. This work was supported by the ERC consolidator grant CUTOFF (101123174). 
\bibliographystyle{plain}
\bibliography{draft}
\end{document}